\documentclass[conference]{IEEEtran}
\usepackage{times}

\let\labelindent\relax
\usepackage[numbers]{natbib}
\usepackage{multicol}
\usepackage{enumitem}
\usepackage{graphicx}
\usepackage{amssymb,mathtools,amsfonts}
\usepackage{gensymb}
  
\usepackage{amsthm}  
\usepackage{caption, subcaption}
\usepackage[bookmarks=true]{hyperref}
\usepackage{cleveref}
\usepackage{multirow}
\usepackage{color}
\usepackage{psfrag}
\usepackage{empheq}
\usepackage{bm} 
\usepackage{url}
\usepackage{float}
\usepackage{comment}
\usepackage{alphalph}
\usepackage[short]{optidef}
\usepackage[dvipsnames]{xcolor}
\usepackage{tikz}
\usepackage{suffix}

\usetikzlibrary{positioning}
\usetikzlibrary{arrows}
\usetikzlibrary{babel}
\usetikzlibrary{automata}
\usetikzlibrary{shapes, snakes}
\usetikzlibrary{positioning, fit, calc}
\usetikzlibrary{decorations.pathreplacing,calligraphy}
\usetikzlibrary{patterns}
\usetikzlibrary{patterns.meta}
\usepackage{pgfplots}
\pgfplotsset{compat=newest}
\usepackage{pgfplotstable}
\usepackage{pgf}

\newcommand{\nosnoc}{\texttt{nosnoc}}

\newcommand{\R}{{\mathbb{R}}}
\newcommand{\dd}{{\mathrm{d}}}

\newcommand{\ns}{{n_\mathrm{s}}}

\newcommand{\Nfe}{{N_\mathrm{fe}}}
\newcommand{\sign}{{\mathrm{sign}}}

\newcommand{\zsdf}{{z_{\mathrm{sdf}}}}
\newcommand{\gsdf}{{g_{\mathrm{sdf}}}}
\newcommand{\normal}{d_n}
\newcommand{\tangent}{d_t}
\newcommand{\lnormal}{\lambda}
\newcommand{\ltangent}{\lambda^t}
\newcommand{\ltangentp}{\lambda^{t+}}
\newcommand{\ltangentn}{\lambda^{t-}}
\newcommand{\Gf}{G^f}
\newcommand{\Hf}{H^f}
\newcommand{\vt}{v^{t}}

\newcommand{\paren}[1]{\left( #1 \right)}
\newcommand{\brac}[1]{\left[ #1 \right]}
\newcommand{\cbrac}[1]{\left\{ #1 \right\}}
\WithSuffix\newcommand\paren*[1]{( #1 )}
\WithSuffix\newcommand\brac*[1]{[ #1 ]}
\WithSuffix\newcommand\cbrac*[1]{\{ #1 \}}

\newcommand{\vo}{\vec{o}\@ifnextchar{^}{\,}{}}

\newcommand{\Set}[2]{\left\{\, #1 \mid #2\,\right\}}
\WithSuffix\newcommand\Set*[2]{\{\, #1 \mid #2\,\}}
\newcommand{\norm}[1]{{\left\Vert#1\right\Vert}}
\newcommand{\transp}[1]{{#1}^\top}

\newcommand{\C}{{\mathcal{C}}}
\newcommand{\xdot}{{\dot{x}}}
\newcommand{\proj}[2]{\mathrm{P}_{#1}\paren{#2}}
\WithSuffix\newcommand\proj*[2]{\mathrm{P}_{#1}\paren*{#2}}
\newcommand{\tancone}[2]{\mathcal{T}_{#1}\paren{#2}}

\DeclareMathOperator*{\diag}{diag}
\DeclareMathOperator*{\argmin}{arg\,min}

\begin{document}

\title{Towards Solutions of Manipulation Tasks via Optimal Control of Projected Dynamical Systems}

\author{\authorblockN{Anton Pozharskiy, Armin Nurkanovi\'c, Moritz Diehl}
  \authorblockA{University of Freiburg\\
    Freiburg, Germany 79110\\
  firstname.lastname@imtek.uni-freiburg.de\vspace{-0.5cm}}}

\maketitle
\begin{abstract}
  We introduce a modeling framework for manipulation planning based on the formulation of the dynamics as a projected dynamical system.
  This method uses implicit signed distance functions and their gradients to formulate an equivalent gradient complementarity system.
  The optimal control problem is then solved via a direct method, discretized using finite-elements with switch detection\cite{Pozharskiy2024}.
  An extension to this approach is provided in the form of a friction formulation commonly used in quasi-static models.
  We show that this approach is able to generate trajectories for problems including multiple pushers, friction, and non-convex objects modeled as unions of convex ellipsoids with reasonable computational effort.
\end{abstract}

\IEEEpeerreviewmaketitle
\section{Introduction}
Optimal control of systems with contacts is a difficult problem in the field of robotic manipulation.
This is due to the fact that these kinds of systems are fundamentally under-actuated and encompass nonsmooth phenomena such as contact forces and friction.
The under-actuation arises from the fact that the unactuated slider cannot be arbitrarily moved in a given state of the system.
The nonsmooth, hybrid dynamics occur due to the switches from contact constrained motion to contact-free motion as well as the sign of friction forces.

Much work has been done in this field using approaches that mix of discrete and continuous elements.
This includes work using graphs of convex sets~\cite{Graesdal2024}, rapidly exploring random trees~\cite{Zito2012}, and differential dynamic programming~\cite{Xue2023}.

Projected Dynamical Systems (PDS) and a related class of systems called First Order Sweeping Processes (FOSwP) have had relatively little study in the context of optimal control.
Some theoretical work has been done on the optimal control of FOSwP with non-overlapping disc conditions~\cite{Cao2021,Cao2021a} and some numerical approximations for more generic problems~\cite{Pinho2020}.
Specifically for PDS there is little optimal control literature however they have used for control synthesis~\cite{Hauswirth2021,Lorenzetti2022} and were originally proposed in the study of economics~\cite{Dong1996,Nagurney1995}.

We propose a framework which tackles manipulation problems using a PDS to model its dynamics.
We then apply direct optimal control to an equivalent Dynamical Complementarity System (DCS) by first accurately discretizing the Optimal Control Problem (OCP) then solving the discrete-time OCP.
The resulting finite dimensional optimization problem is a mathematical program with complementarity constraints which is solved via a relaxation homotopy method~\cite{Nurkanovic2024b}.

\section{Projected dynamical systems}
In this section we briefly review the mathematical formulation we will use to model the manipulation tasks in later sections.
First we define the tangent cone to $\C$ at $x$, denoted by $\tancone{\C}{x}$, as the set of all vectors $d\in \R^n$ for which there exists sequences $\{x_i\} \subset \C$ and $\{t_i\}, t_i \ge 0$, with $x_i \to x$ and $t_i\to0$, such that ${d = \lim_{i\to \infty} \frac{x_i-x}{t_i}}$.
We also define the projection operator: ${\proj{K}{x} =\argmin_{s \in K}\frac{1}{2}\norm{s-x}^2_2}$.
Using these preliminaries we define a PDS which is a form of constrained dynamical system where the state is limited to evolving within a set.
The dynamics of such a system are written as:
\begin{equation}
  \label{eq:PDS}
  \xdot(t) = \proj{\tancone{\C}{x(t)}}{f(x(t))},
\end{equation}
with $x(0) \in \C$ and $f:\R^{n_x}\rightarrow\R^{n_x}$ which is at least twice continuously differentiable.
For computational tractability we treat sets written as $\C = \Set*{x\in\R^n}{c(x)\ge 0}$ with ${c:\R^{n_x}\rightarrow \R^{n_c}}$.
In order to numerically solve the PDS we formulate it as an equivalent DCS~\cite{Brogliato2006}:
\begin{subequations}
  \begin{align*}
    \xdot(t) &= f(x(t)) + \nabla c(x(t))\lambda(t),\\
    0 &\le c_i(x(t)) \perp \lambda_i(t) \ge 0, \quad i=1,\ldots,n_c. \label{eq:GCS:comps}
  \end{align*}
\end{subequations}
In the manipulation setting the functions $c(x(t))$ can be interpreted as signed distance functions between objects and the multipliers $\lambda$ as contact forces which act to prevent object overlap.

\section{Modeling manipulation tasks}
In this section we describe the method for modeling planar manipulation tasks using spherical approximations.
This is followed with a promising direction for modeling manipulation tasks with ellipsoidal approximations.
Finally we discuss an extension of the simple PDS equivalent model with the well studied quasi-static friction model~\cite{Graesdal2024, Halm2020}.

\subsection{Spherical manipulation tasks}
We model the sliding contact dynamics between two discs in the plane via a non-overlapping constraint.
In this case we treat two circles with centers $p_1,p_2 \in \R^2$ and radiuses $r_1$ and $r_2$, respectively.
The non-overlapping of these two discs is written as the inequality $\norm{p_1-p_2}^2-(r_1+r_2)^2 \ge 0$.
The left hand side expression can be interpreted as a signed distance function between the two discs which must be kept non-negative as in the finite definition of the set $\C$.
We model any number of these pairwise contacts via these signed distance functions $c_{ij}(x) = \norm{p_i-p_j}^2-(r_i+r_j)^2$.
In the case of frictionless sliding contacts the existence and uniqueness of the solutions to such a system can be guaranteed if the gradients of the constraint function $c(x)$ satisfy the linear-independence constraint qualification, i.e. the columns of $\nabla_xc(x)$ are linearly independent at all points $x$.

\subsection{Ellipsoidal manipulation tasks}
Moving from the simplicity of the spherical manipulation tasks in the prequel to using ellipsoids to model objects requires new machinery to define the signed distance distance which we can embed within our formulation.
There exist disc-ellipsoid and ellipsoid-ellipsoid distance functions that can be calculated as the solutions to polynomials of various orders though often one needs to resolve which solution is the correct distance~\cite{Uteshev2018}.
As such, we choose a formulation in the form of a convex program which has a unique solution and a unique contact point.

First, we define an ellipse augmented with its center $p\in\R^2$ and rotation $\theta\in\R$ using the set $E(p,\theta) = \Set*{x\in \R^2}{\transp{(x-p)}P(\theta)(x-p) \le 1},$ where $P(\theta)=\transp{R(\theta)}\hat{P}R(\theta)$, $R(\theta)$ is the 2d rotation matrix and $\hat{P}\in\R^{2\times 2}$ is a positive definite symmetric matrix.
Using this definition we devise a signed distance function between two ellipses which is the solution to the convex optimization problem~\cite{Tracy2023}:
\begin{mini*}[3]
  {\substack{\alpha\in\R,\\p_d\in\R^2}}
  {\alpha -1}{\label{eq:ell_disc}}
  {c(x) =}
  \addConstraint{\alpha}{\le\transp{(p_d-p_i)}P(\theta_i)(p_d-p_i),\ }{i=1,2,}
\end{mini*}
where $x= (p_1,\theta_1,p_2,\theta_2)\in\R^6$.
In this optimization problem, $\alpha$ represents a factor by which each ellipse is rescaled.
The distance is found by minimizing $\alpha$ such that there is still a common point, the contact point $p_d$, between the two rescaled ellipses, enforced by the constraints now collected as $g(\alpha,p_d;x)$.
We subtract one in the objective so that if $\alpha \le 1$, $c(x) \le 0$ and there is an overlap between the original objects.
Due to the convexity of this problem we embed it using its KKT conditions:
\begin{align*}
  1-\mu_1-\mu_2 &= 0,\\
  2\mu_1\transp{P(\theta_1)}(p_d-p_1)+ 2\mu_2\transp{P(\theta_2)}(p_d-p_2) &=0,\\
  0 \le \mu_i \perp e\transp{(p_d-p_i)}P(\theta_i)(p_d-p_i) - \alpha&\ge 0,\ i=1,2,
\end{align*}
introducing the Lagrange multipliers $\mu\in\R^2$.
In order to use this formulation of the signed distance function (SDF), we extract from it its gradient $\normal(x) = \nabla_xc(x) = \nabla_x(g(\alpha,p_d;x))\mu$ as in \cite{Dietz2024}.
The contact point is unique and therefore for any configuration of two ellipsoids $\normal(x)$ is unique.
This is sufficient to show that the tangent cone must be closed and convex at all points~\cite{Rockafellar2004}.
As such, we can show that unique solutions exist (under mild assumptions on $f(x,u)$) for the projected dynamical system in \Cref{eq:PDS}~\cite{Nagurney1995}.

We can further model non-convex shapes using the union of multiple ellipses with a common point used as the center of each ellipse.
It remains to be shown whether this formulation maintains the required properties of the set $\C$.
However, experience with a variety of problems suggests that these systems are well posed, including in the case of multiple contacts in interior corners of unions-of-ellipsoids objects.

\subsection{Modeling friction}\label{sec:friction}
In order to model friction we depart the standard theoretical model of projected dynamical systems.
We apply a common approach by modeling friction via a quasi-static model~\cite{Trinkle1989,Xue2023}, i.e., we expect the friction force to minimize the relative tangential velocity between contacting objects.
To this end $\ltangentp$, $\ltangentn$, and $\gamma$ are introduced as the positive and negative components of the friction force and an additional slack to model the maximum friction force respectively.
In the planar case we can calculate the tangential velocity by projecting the velocity onto the tangent line $\tangent = R(\pi)\normal$.
For brevity we introduce a variable ${\vt_{i} = \frac{\dot{p}_i\cdot\tangent}{\tangent\cdot\tangent} + \dot{\theta}_i\norm{p_d - p_i}}$; the tangential velocity of the contact point for each object in each contact.
The quasi-static friction formulation requires the following conditions: ${\ltangent=\ltangentp-\ltangentn}$ and ${\ltangent = -\mu_f\lnormal\sign(\vt_2-\vt_1)}$.
We augment the system in \Cref{eq:PDS} with the following equivalent complementarity conditions:
\begin{subequations}
  \begin{eqnarray}
  \label{eq:friction_comp}
    0 \le \vt_{1} - \vt_{2} + \gamma \perp \ltangentp \ge 0,\\
    0 \le \vt_{2} - \vt_{1} + \gamma \perp \ltangentn \ge 0,\\
    0 \le \mu_f\lnormal - \ltangentp - \ltangentn \perp \gamma \ge 0.
  \end{eqnarray}
\end{subequations}
$z = (\lnormal,\ltangentp,\ltangentn,\gamma,\vt)$ summarizes all the algebraic variables.
This yields the expression for the dynamics ${\hat{f}(x, u, z) = f(x,u) + \normal\lnormal - \tangent\ltangentp + \tangent\ltangentn}$.
For compactness of notation, complementarity pairs for all friction components in \Cref{eq:friction_comp} are collated into ${0\le\Gf\perp\Hf\ge 0}$.
\section{Solving the OCP}
\begin{figure*}[t]
  \centering
  \begin{subfigure}{0.19\linewidth}
    \includegraphics[width=\linewidth]{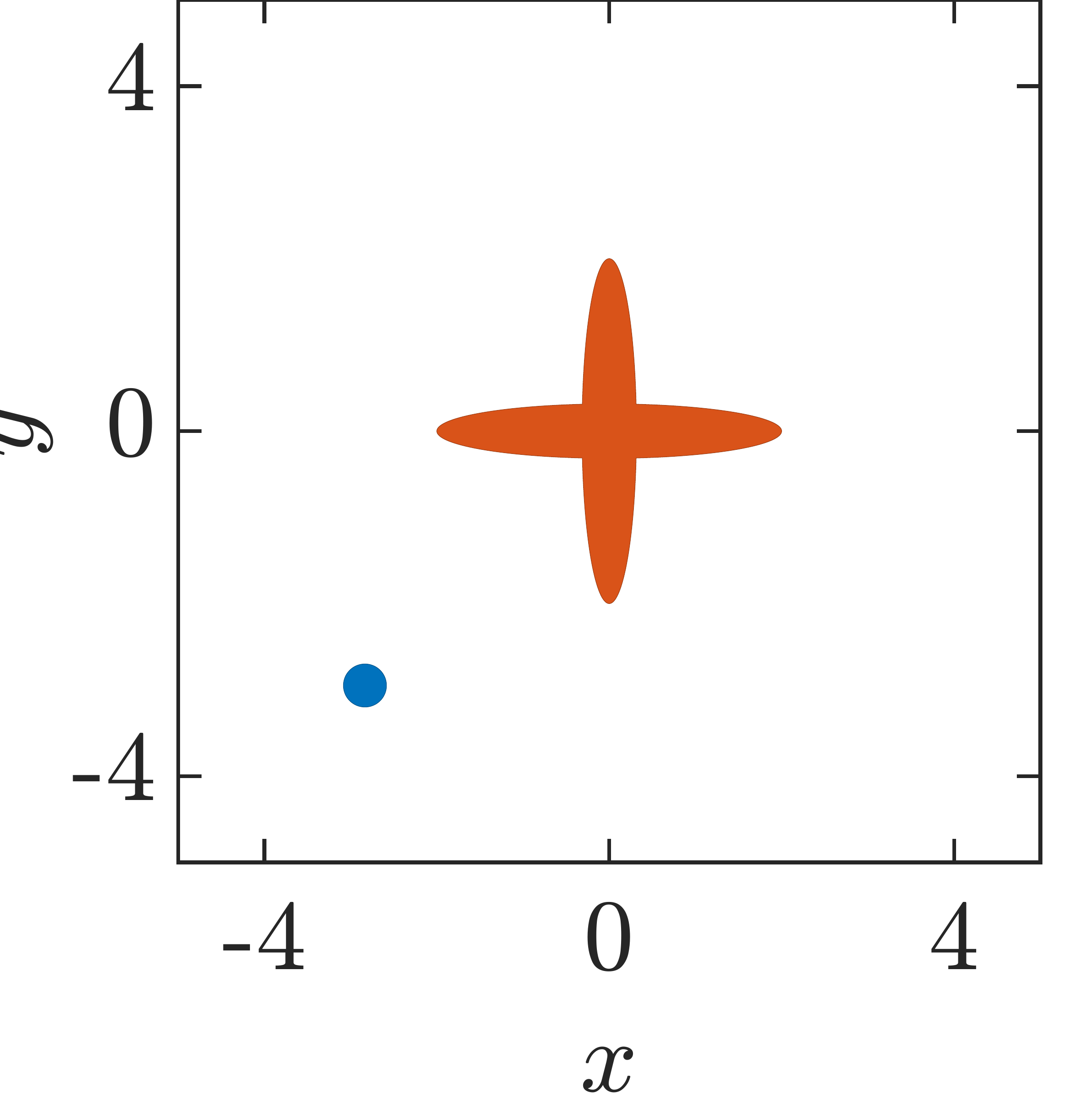}
    \caption{$t=0$}
    \label{fig:spin:1}
  \end{subfigure}
  \begin{subfigure}{0.19\linewidth}
    \includegraphics[width=\linewidth]{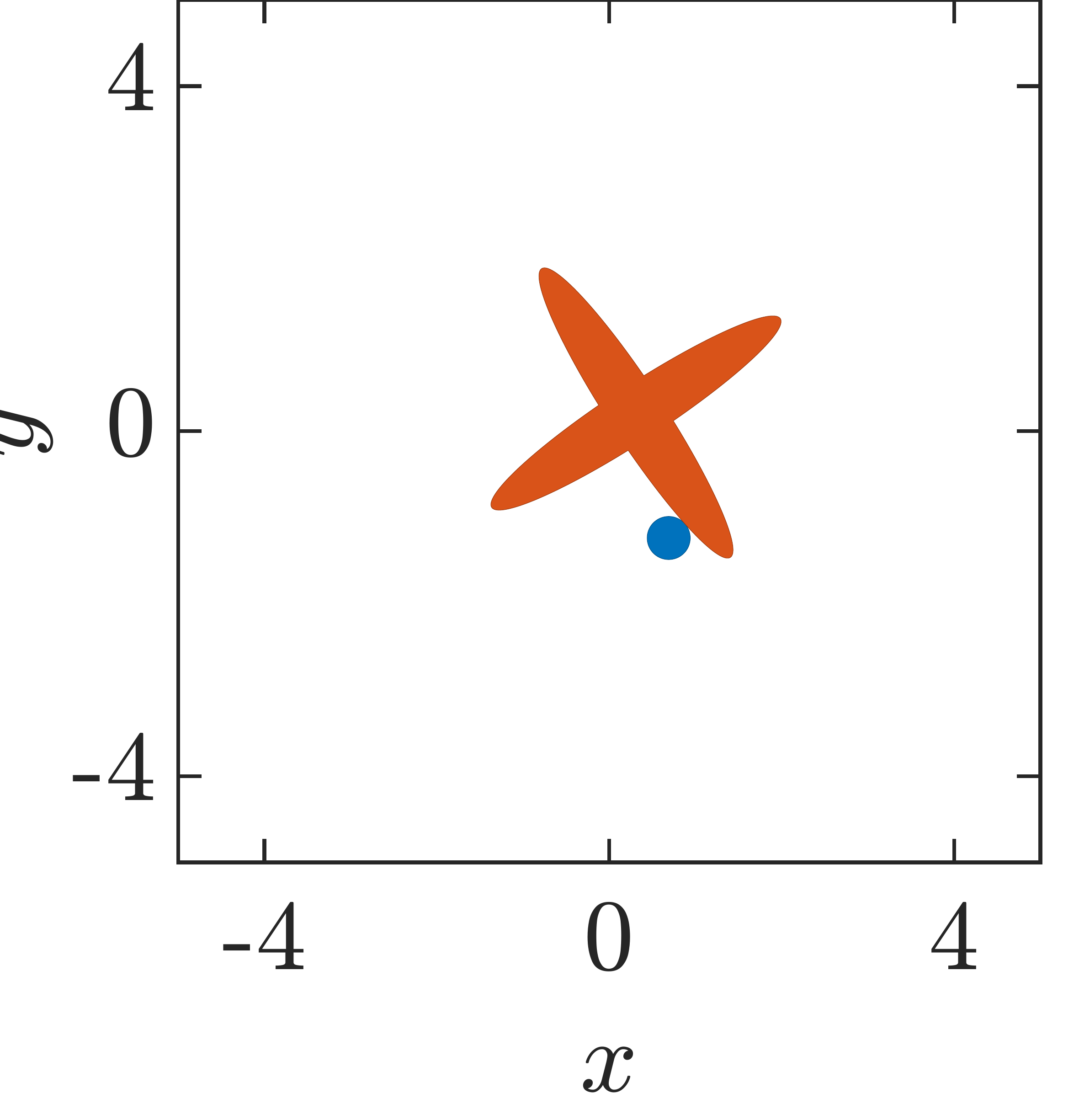}
    \caption{$t=3.833$}
  \end{subfigure}
  \begin{subfigure}{0.19\linewidth}
    \includegraphics[width=\linewidth]{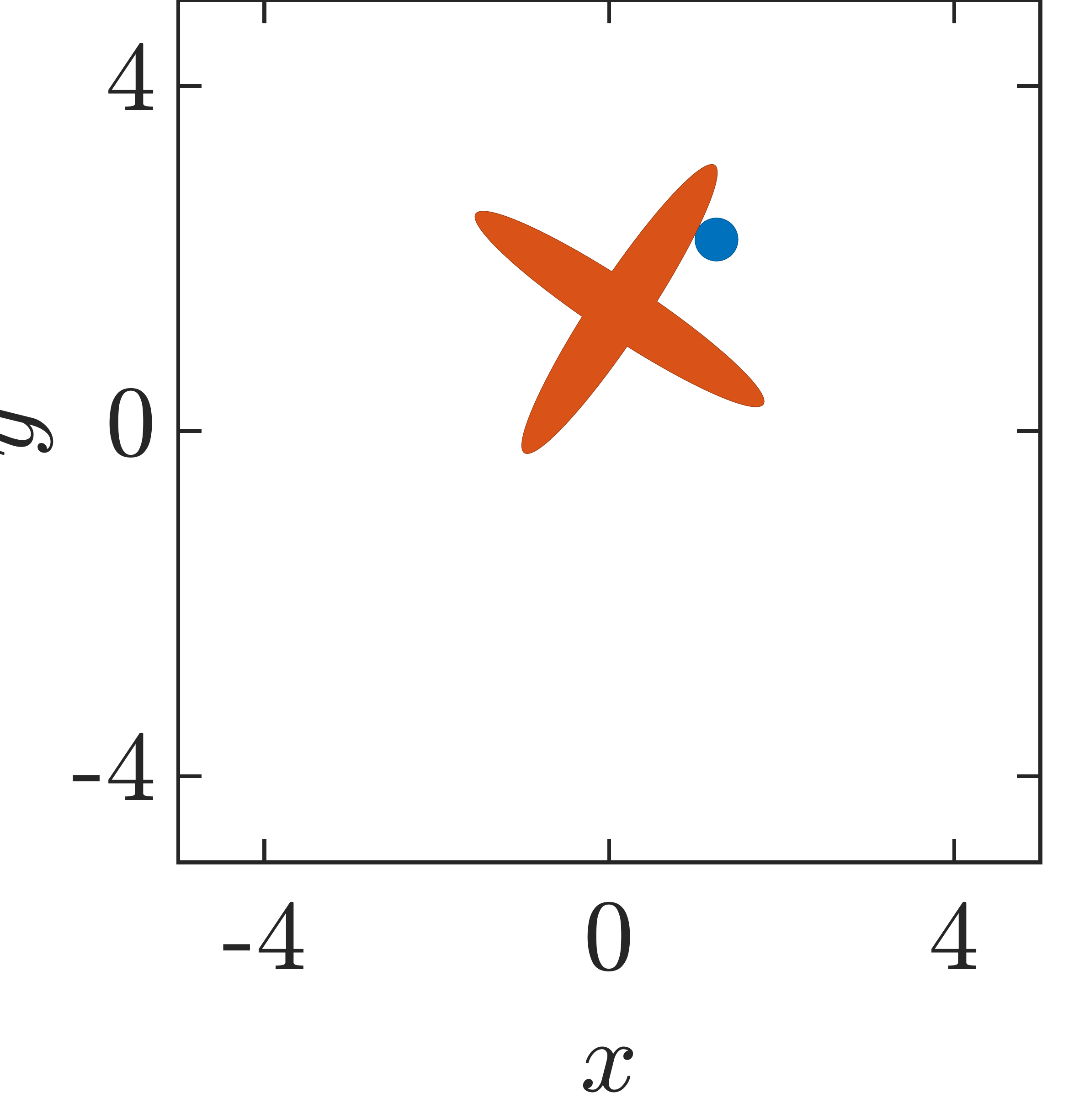}
    \caption{$t=8.5$}
  \end{subfigure}
  \begin{subfigure}{0.19\linewidth}
    \includegraphics[width=\linewidth]{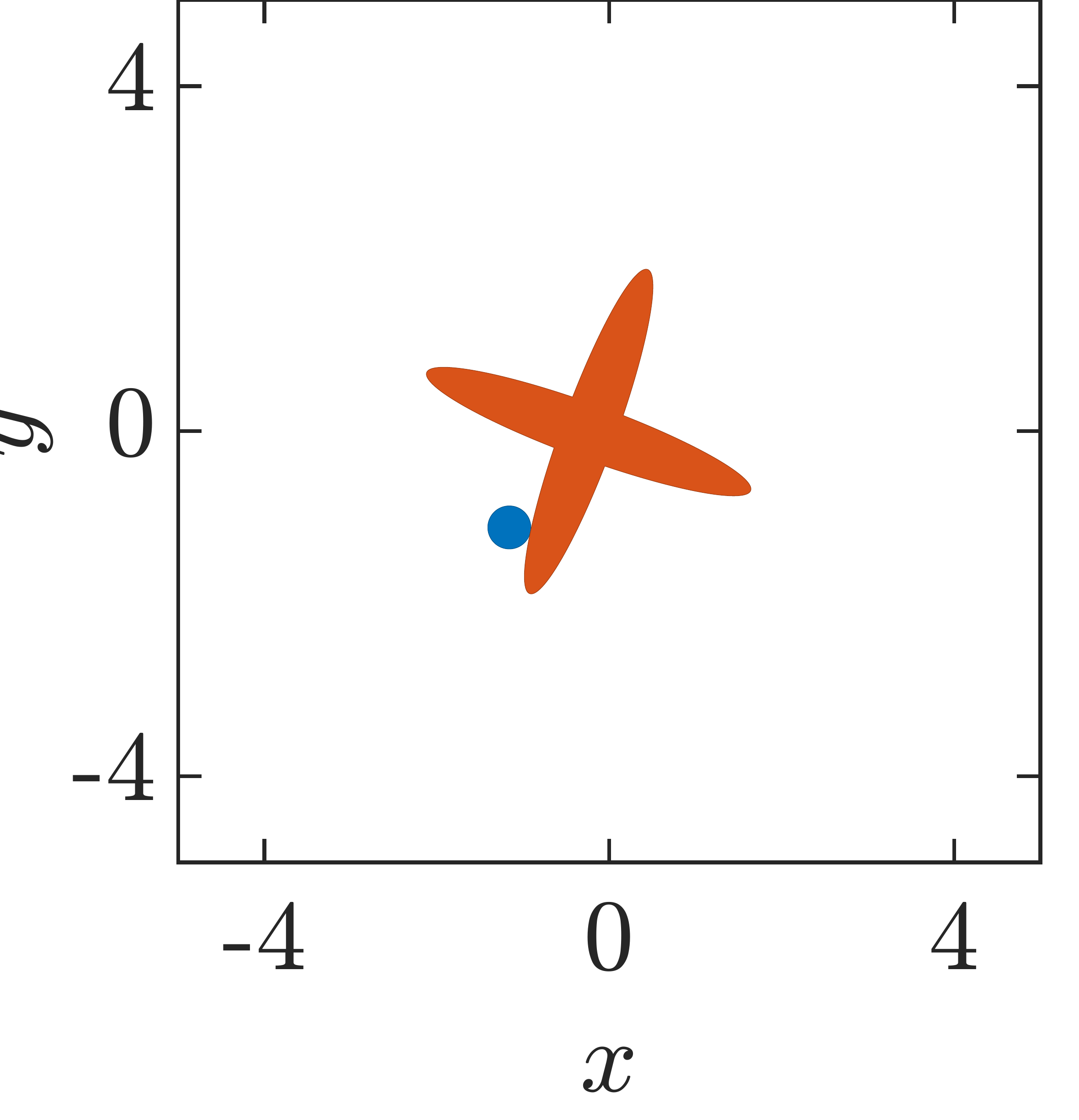}
    \caption{$t=16.667$}
  \end{subfigure}
  \begin{subfigure}{0.19\linewidth}
    \includegraphics[width=\linewidth]{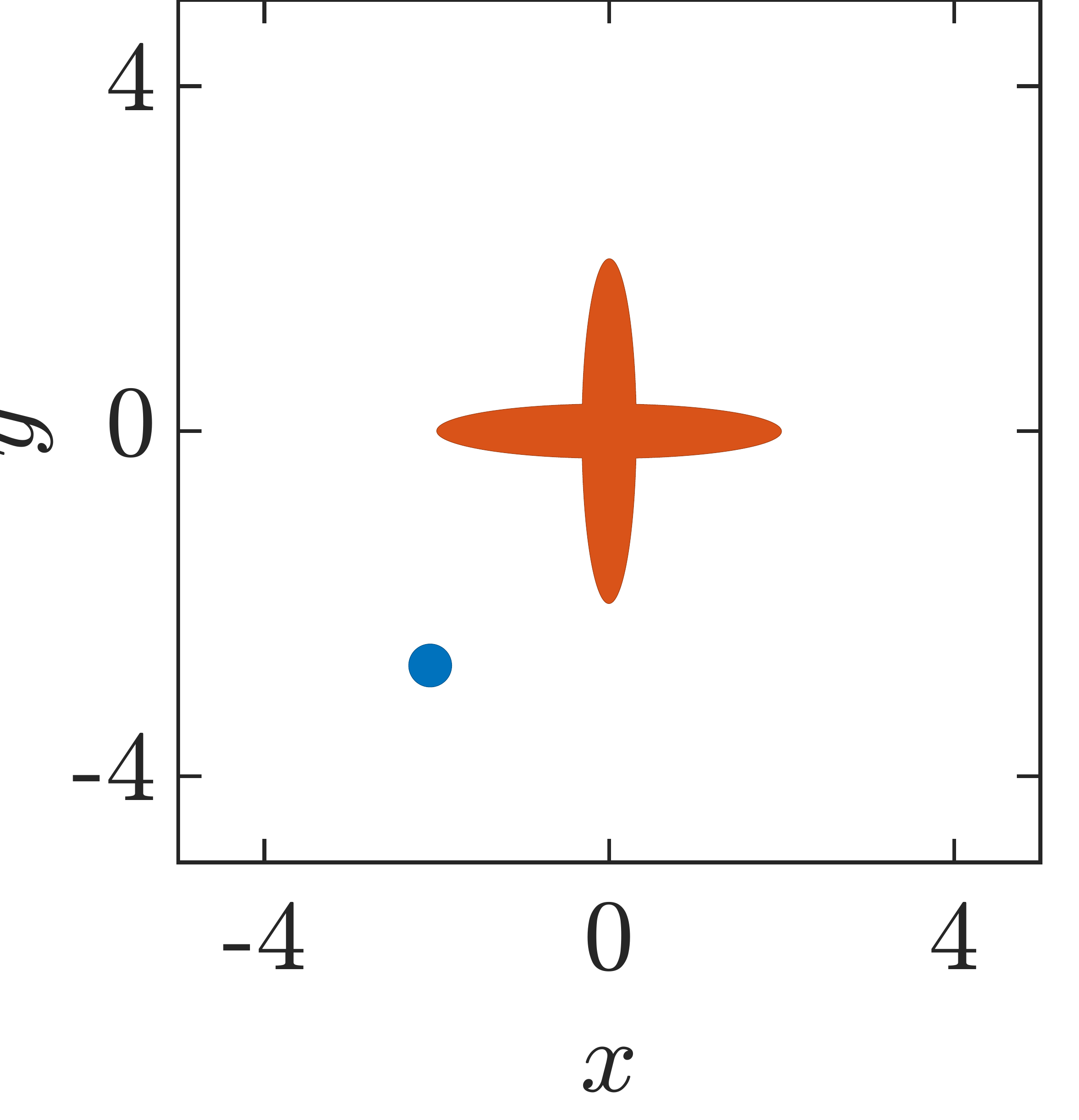}
    \caption{$t=20.0$}
  \end{subfigure}
  \vspace{-0.2cm}
  \caption{Several frames of the solution for the frictionless manipulation problem.}
  \vspace{-0.7cm}
  \label{fig:spin}
\end{figure*}
\begin{figure}[t]
  \centering
  
  \includegraphics[width=\linewidth]{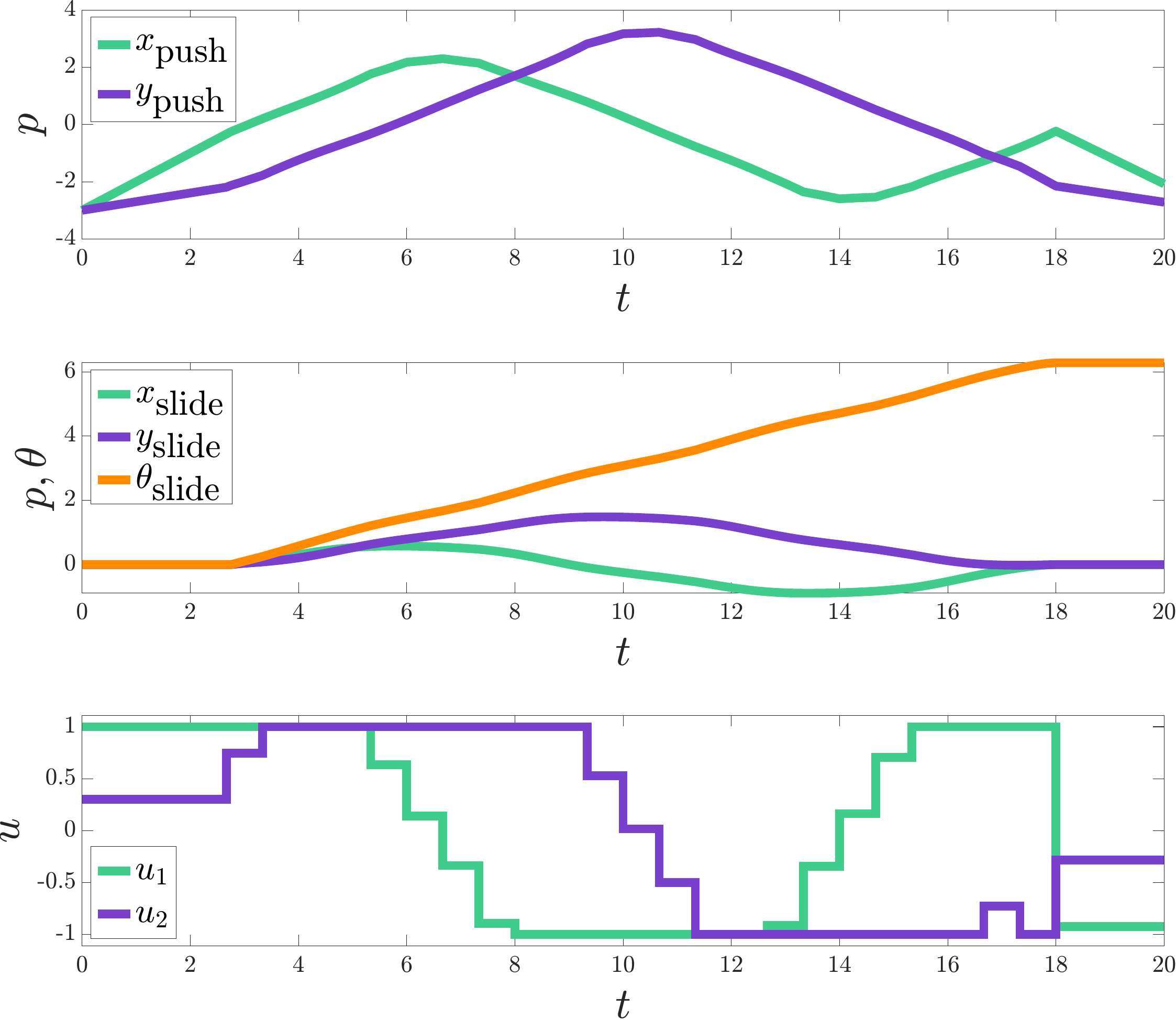}
  
  \caption{Plots of solution to the frictionless pushing problem.}
  \vspace{-0.7cm}
  \label{fig:spinning-states}
\end{figure}
In this section we discuss the approach we take to accurately solving the Optimal Control Problem (OCP) that results from the above models, via direct optimal control.
\subsection{Finite Elements with Switch Detection}
\begin{figure*}[t]
  \centering
  \begin{subfigure}{0.19\linewidth}
    \includegraphics[width=\linewidth]{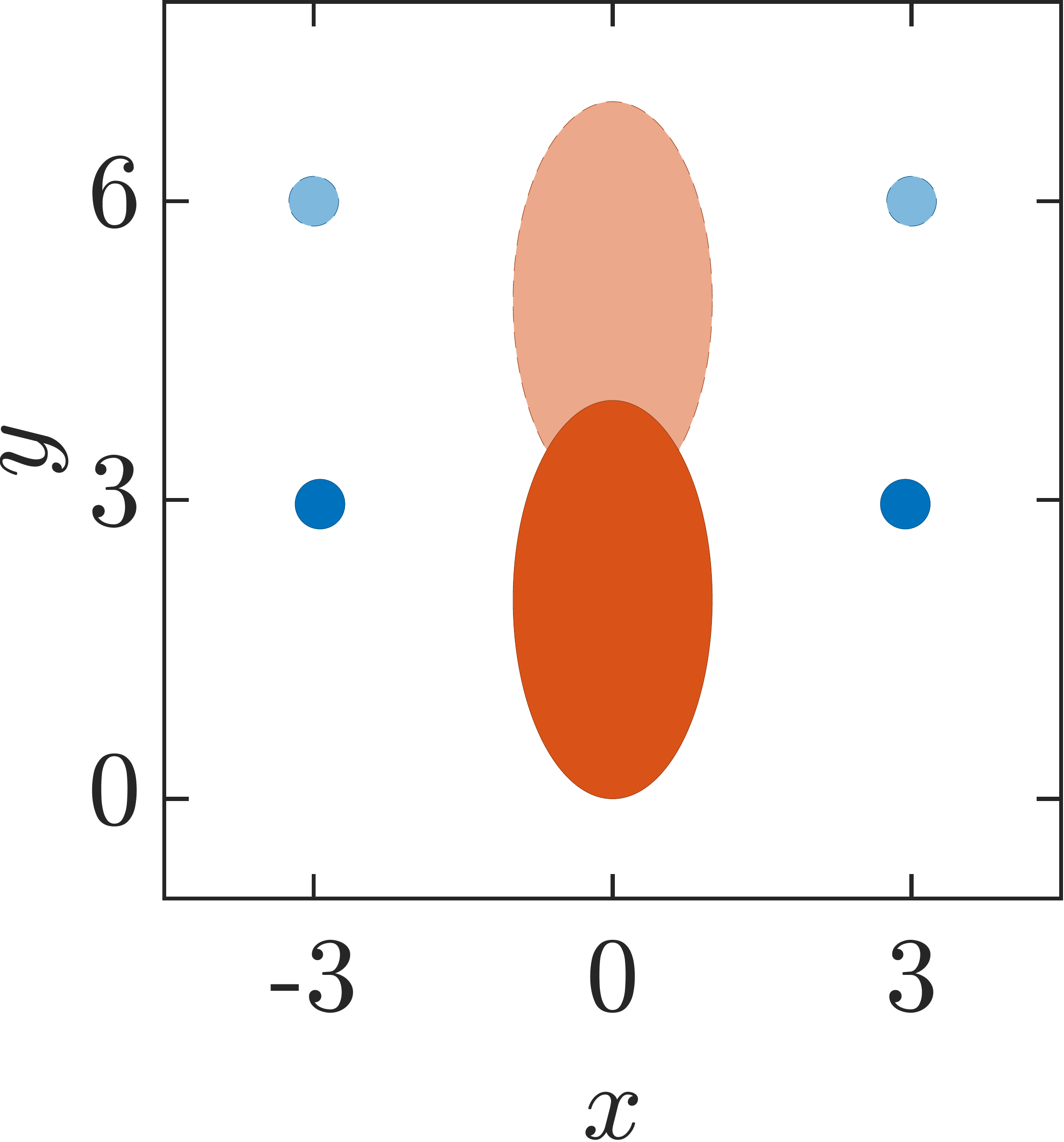}
    \vspace{-0.4cm}
    \caption{$t=0.0$}
    \label{fig:drag:1}
  \end{subfigure}
  \begin{subfigure}{0.19\linewidth}
    \includegraphics[width=\linewidth]{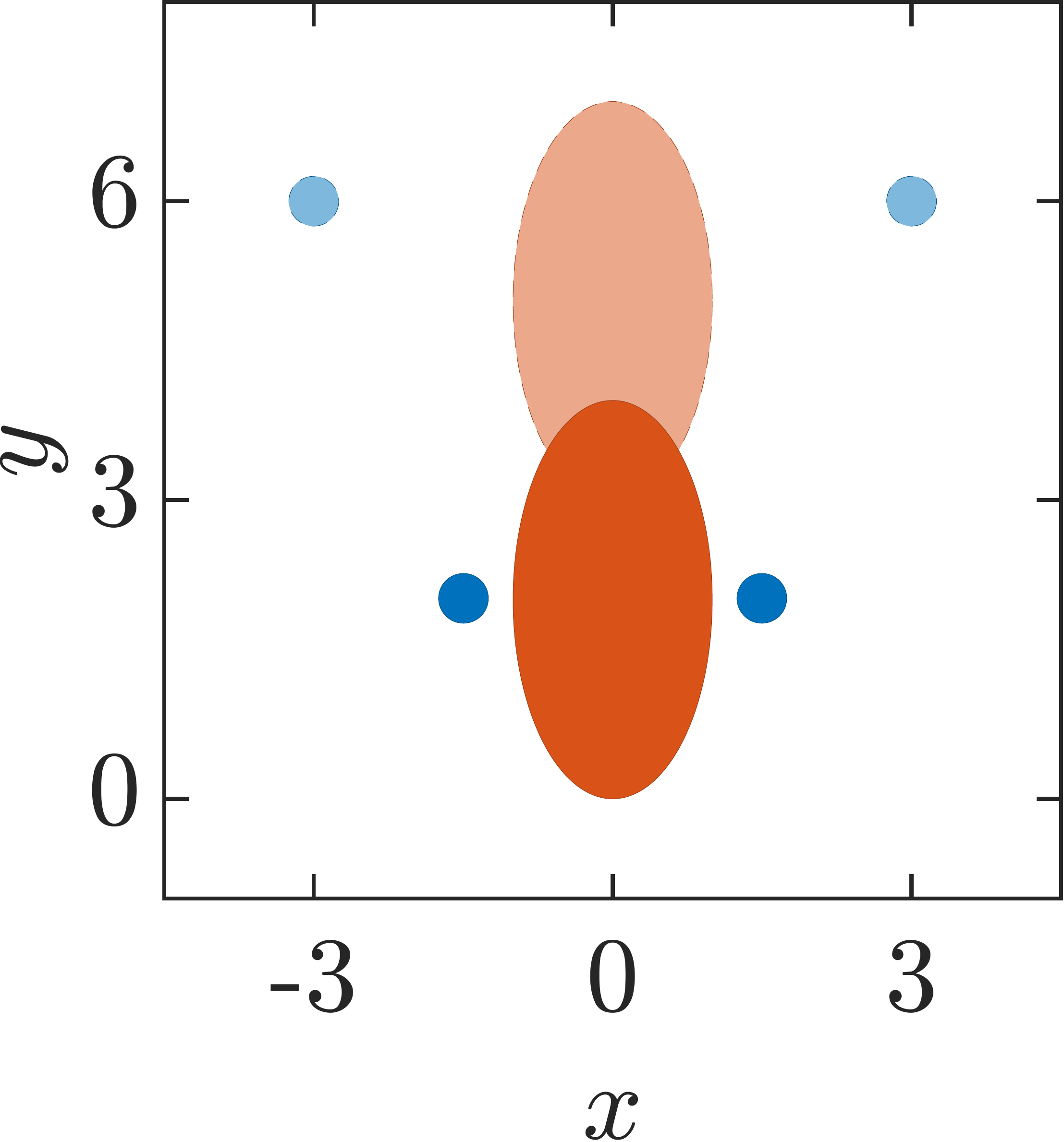}
    \vspace{-0.4cm}
    \caption{$t=3.833$}
  \end{subfigure}
  \begin{subfigure}{0.19\linewidth}
    \includegraphics[width=\linewidth]{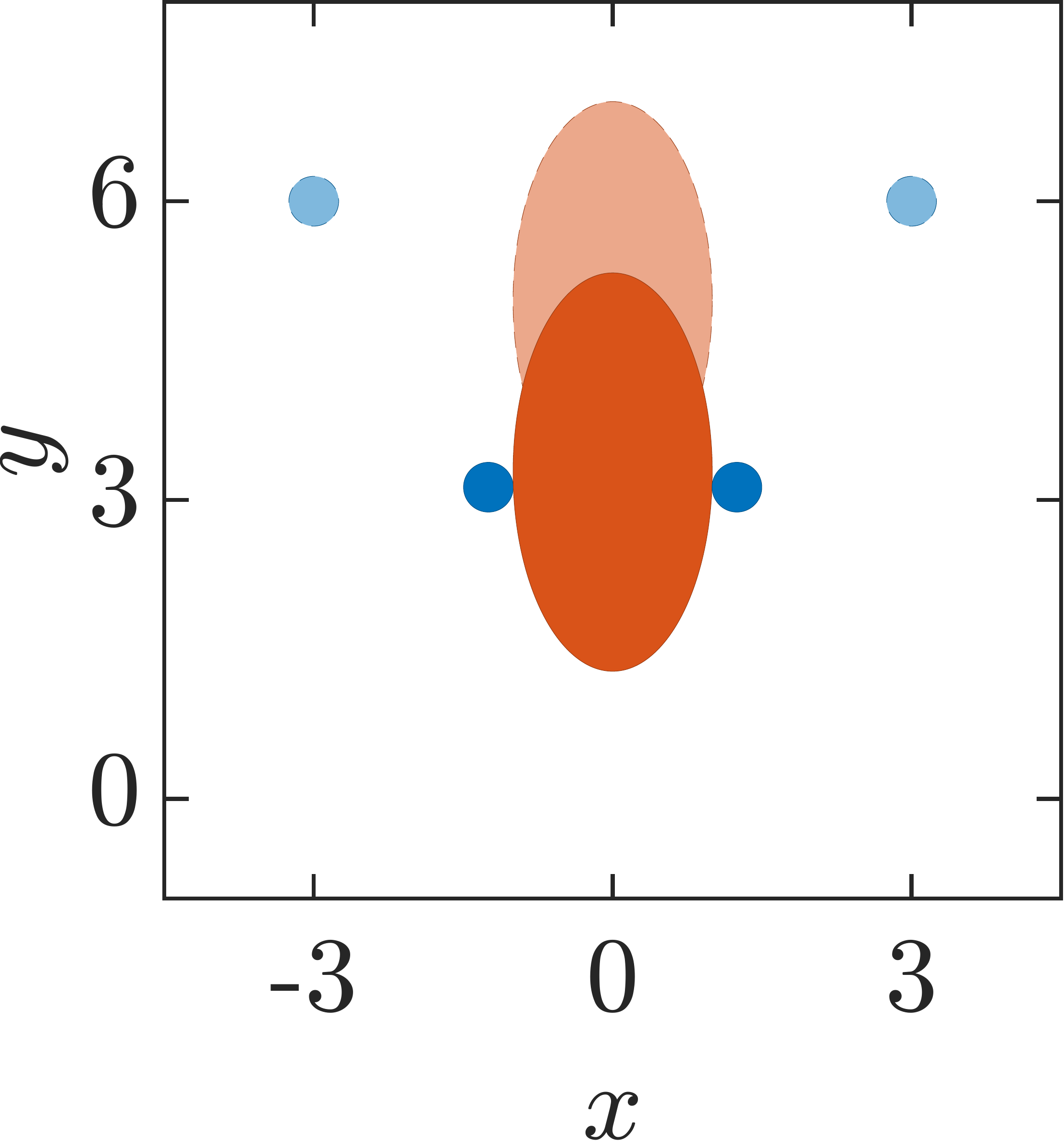}
    \vspace{-0.4cm}
    \caption{$t=8.5$}
  \end{subfigure}
  \begin{subfigure}{0.19\linewidth}
    \includegraphics[width=\linewidth]{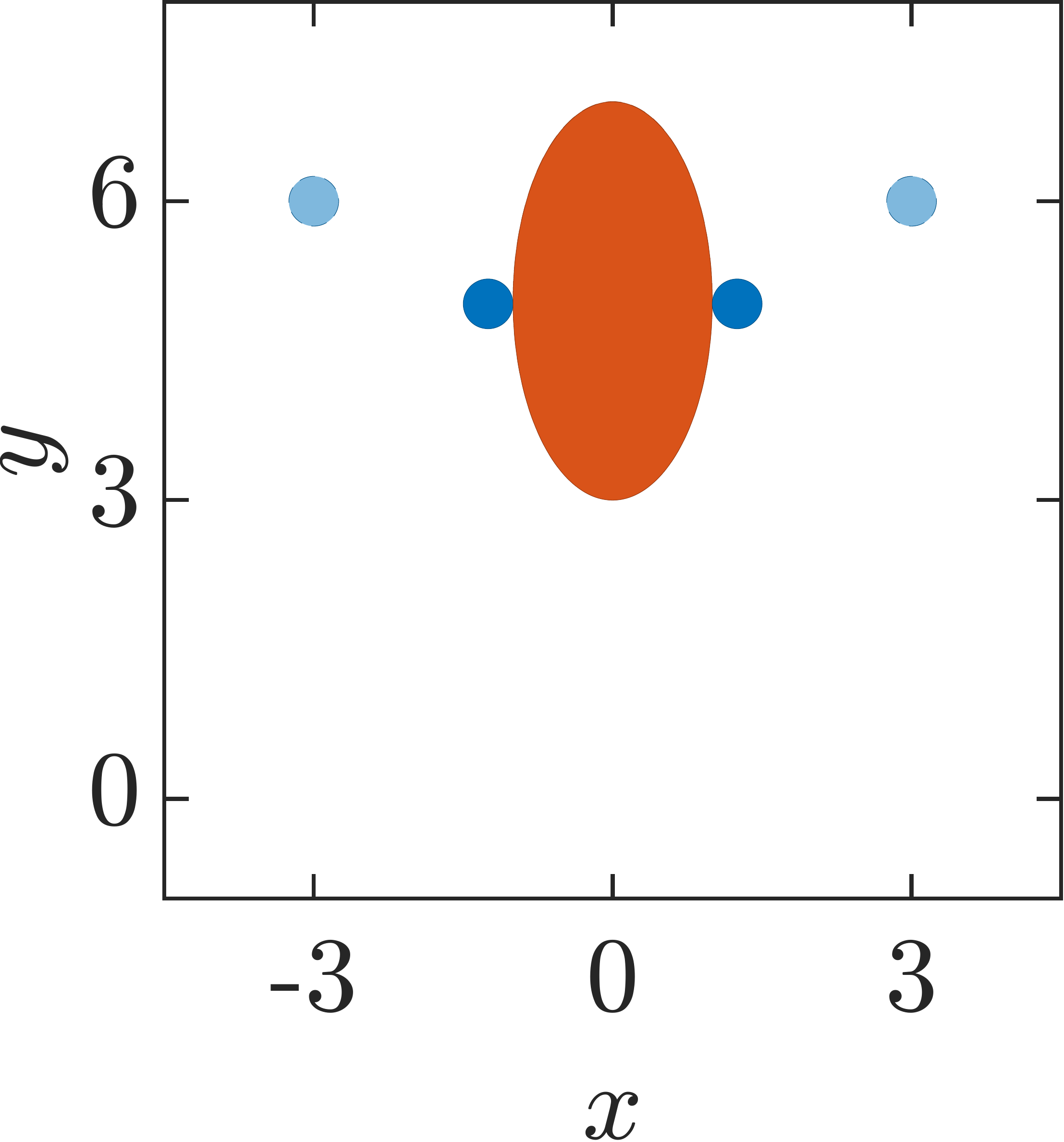}
    \vspace{-0.4cm}
    \caption{$t=14.5017$}
  \end{subfigure}
  \begin{subfigure}{0.19\linewidth}
    \includegraphics[width=\linewidth]{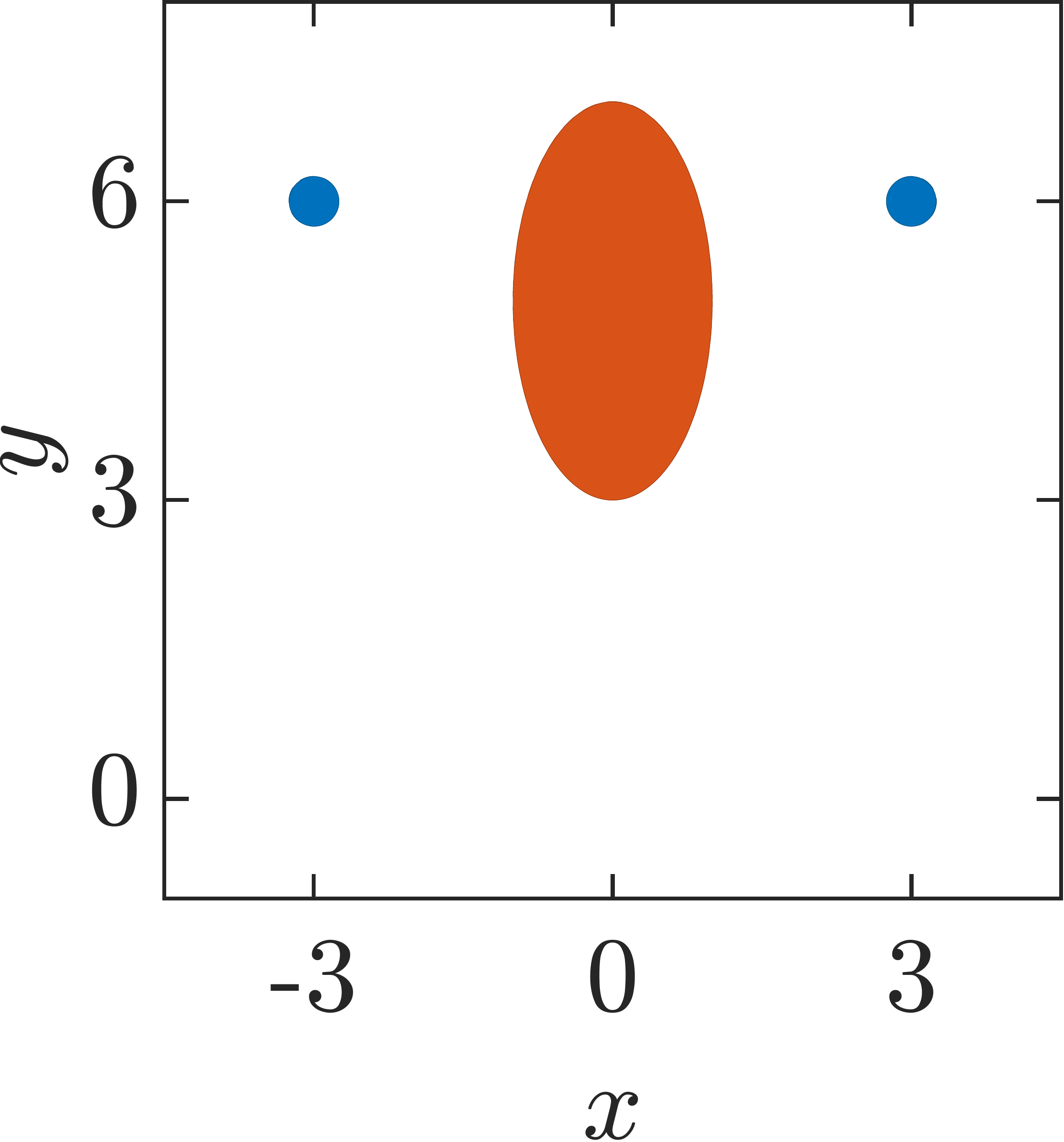}
    \vspace{-0.4cm}
    \caption{$t=20.0$}
  \end{subfigure}
  \vspace{-0.2cm}
  \caption{Several frames of the solution for the manipulation problem with $\mu_f=0.5$.}
  \vspace{-0.8cm}
  \label{fig:drag}
\end{figure*}
We apply the Finite Elements with Switch Detection~(FESD) discretization for Projected Dynamical Systems~(PDS)~\cite{Pozharskiy2024} to discretize an OCP in the form:
\begin{mini*}[3]
  {x(\cdot),u(\cdot), z(\cdot),\zsdf(\cdot)}
  {\int_{0}^{T} L(x(t),u(t))\dd t +  M(x(T))}{\label{eq:ocp}}
  {}
  \addConstraint{x(0)}{=x_0}
  \addConstraint{\xdot(t)}{= \hat{f}(x(t), u(t), z(t))),\ }{t \in [0, T]}
  \addConstraint{0}{\leq c(x(t))\perp \lnormal(t)\geq 0,\ }{t \in [0,T]}
  \addConstraint{0}{\leq \Gf(t)\perp \Hf(t)\geq 0,\ }{t \in [0,T]}
  \addConstraint{0}{\le \gsdf(x(t),\zsdf(t)),\ }{t \in [0,T]}
  \addConstraint{0}{\geq g_{\mathrm{p}}(x(t),u(t)),\ }{t \in [0,T]}
  \addConstraint{0}{\geq g_{\mathrm{t}}(x(T)).}
\end{mini*}
With the stage cost term $L:\R^{n_x}\times\R^{n_u}\rightarrow\R$, terminal cost term $M:\R^{n_x}\rightarrow\R$, initial state ${x_0\in\R^{n_x}}$, path constraints $g_{\mathrm{p}}:\R^{n_x}\times\R^{n_u}\rightarrow\R^{n_{g_\mathrm{p}}}$, friction complementarities $\Gf,\ \Hf\in \R^{5n_c}$, terminal constraints ${g_{\mathrm{t}}:\R^{n_x}\rightarrow\R^{n_{g_\mathrm{t}}}}$, controls $u(t)\in\R^{n_u}$, and SDF algebraics ${\zsdf(t) = (\alpha, p_d)\in\R^{3n_c}}$.
We discretize the each of $N_s$ control intervals via $\Nfe$ finite elements which is each a Runge-Kutta discretization with $\ns$ stages and a Butcher tableau~\cite{Hairer1996} with entries $a_{i,j}$, $b_j$ for $i,j = 1,\ldots,\ns$.
\begin{align*}
  &x_{n+1,0} = x_{n,0} + h_n\sum_{j=1}^{\ns}b_j\hat{f}(x_{n,j}, u, z_{n,j}),\\
  &x_{n,i}\! = \! x_{n,0} + h_n\sum_{j=1}^{\ns}a_{i,j}\hat{f}(x_{n,j}, u, z_{n,j}),\\
  &0 \le c(x_{n,i}) \perp \lambda_{n,i} \ge 0,\\
  &0 \le \vt_{1,n,i} - \vt_{2,n,i} + \gamma_{n,i} \perp \ltangentp_{n,i} \ge 0,\\
  &0 \le \vt_{2,n,i} - \vt_{1,n,i} + \gamma_{n,i} \perp \ltangentn_{n,i} \ge 0,\\
  &0 \le \mu_f\lnormal_{n,i} - \ltangentp_{n,i} - \ltangentn_{n,i} \perp \gamma_{n,i} \ge 0,
\end{align*}
where $h_n$ are degrees of freedom ${n = 1,\ldots,\Nfe}$, ${i = 1,\ldots,\ns}$, state variables ${x_{n,i} \in\R^{n_x}}$, and algebraics ${z_{i,j}\in \R^{n_c}}$,
Additionally we introduce cross-complementarity constraints which implicitly identify switches in the projected dynamical system:
\begin{align*}
  0 \le c(x_{n,i}) \perp \lambda_{n,j} \ge 0,\\
  0 \le c(x_{n-1,\ns}) \perp \lambda_{n,j} \ge 0,\\
  0 \le c(x_{n,i}) \perp \lambda_{n-1,\ns} \ge 0,
\end{align*}
with $n = 2,\ldots,\Nfe$, $i = 1,\ldots,\ns$, and $j = 1,\ldots,\ns$.

When we consider friction we must also identify switches arising from the system going from stick to slip and vice versa.
Applying the similar approach of complementing the discretization of \Cref{eq:friction_comp} of a finite element with the end of the previous finite element fails.
This is due to the fact that in the case of a contact closing we may have that, for example, both $\gamma$ and $\mu_f\lnormal_{n,i} - \ltangentp_{n,i} - \ltangentn_{n,i}$ are discontinuous.

However, we note that in the absence of a closing contact the functions entering the complementarity conditions in \Cref{eq:friction_comp} are all continuous.
In order to relax the cross complementarity constraints when a contact is closed we take advantage of the fact that a contact is closed at the end of finite element $i$ only if there is a contact $j$ such that $\lambda_{i,j} = c_j(x_i) = 0$.
As such, we introduce a quantity: $b_i = \min(\lambda_i+c(x_i))\ge 0$, for $i=1,\ldots,\Nfe-1$.
The minimum element in the vector $\lambda_i+c(x_i)$ can be found via the linear program:
\begin{mini*}
  {b\in\R}
  {b}{\label{eq:Bmin-opt}}
  {-b_i= }
  \addConstraint{b}{\ge -(\lambda_{i,j}+c_j(x_i)),\ }{j=1,\ldots,n_c}
\end{mini*}
The KKT conditions of which can be embedded by introducing $\mu^b_i\in \R^{n_c}$, and the following $n_c$ complementarities and one equality constraint:
\begin{subequations}
  \begin{align*}
    &1-\sum_{j=1}^{n_c}\mu^b_{i,j}=0\\
    &0\le b+\lambda_i+c(x_i)\perp \mu^b \ge 0
  \end{align*}
\end{subequations}
Using this new quantity the friction cross-complementarities as:
\begin{subequations}
  \begin{align*}
    0\le\Gf_{n,i}\perp\Hf_{n,j}\ge 0\\
    0\le b_n\Gf_{n-1,\ns}\perp\Hf_{n,j}\ge 0\\
    0\le\Gf_{n,i}\perp b_n\Hf_{n-1,\ns}\ge 0
  \end{align*}
\end{subequations}
with $n = 2,\ldots,\Nfe$, $i = 1,\ldots,\ns$, and $j = 1,\ldots,\ns$.
The multiplication by the factor $b_n$ relaxes the continuity conditions when a contact is either closed or opened, as in these cases $\Gf$ and $\Hf$ may both be discontinuous.

\subsection{Solving the discrete problem}
The type of optimization problem produced by the FESD discretization is called a Mathematical Program with Complementarity Constraints (MPCC).
MPCCs are typically solved via relaxation of the complementarity constraints as the original problem is highly degenerate and often difficult to solve directly.
In particular we use a Scholtes relaxation~\cite{Scholtes2001} with a parameter $\sigma$, where the complementarity $G\perp H$ is replaced with $GH\le \sigma$.
A series of relaxed Nonlinear Programs (NLPs) is then solved in a homotopy loop with a decreasing parameter $\sigma$ and warm starting with the solution to the previous iteration.
This is done until the complementarity residual $GH$ is reduced below a given tolerance.
The NLPs are formulated using CasADi~\cite{Andersson2019} and solved using IPOPT~\cite{Waechter2009}.

\section{Experiments}
We now present several examples of interesting manipulation tasks to which we apply our method.
Source code for each of the experiments can be found in the example set of \nosnoc\footnote{On the branch \texttt{pds\_sdf}}~\cite{Nurkanovic2022b}.
For each problem we use the simulation capabilities of \nosnoc~in order to generate a feasible initial guess with zero control input.
In our experience this is particularly important for convergence when using an implicitly calculated gap function $c(x)$.
All experiments were run on a Ryzen 7 5800x with a 3.8 GHz base clock.

\subsection{Planar Pushing}
In this example we model a planar pushing task with frictionless sliding contacts.
This approach can produce non-physical motion as it ignores the friction between the pusher and slider and assumes that they can slide along each other with no frictional forces.
As such, in this setting certain manipulation approaches are precluded such as gripping via two opposing contacts or jamming against static objects.

The state of this system is represented by ${x = (p_1,p_2,\theta_1)\in\R^5}$, the center of mass positions of the pusher and slider followed by the rotation of the slider.
The goal of the problem is for the pusher to rotate the slider a full rotation and return it to its initial position.
We encode this via a quadratic terminal cost ${M(x(T)) = \transp{(x(T)-\bar{x})}Q_T(x(T)-\bar{x})}$ with target state ${\bar{x} = (-3,-3,0,0,2\pi)}$ and weight matrix ${Q_T=\diag(10^{-3}, 10^{-3},10^2,10^2,10^3)}$.
In order to also minimize control effort we include the stage cost term $L(x,u) = \transp{u}Q_uu$ with $Q_u = 10^{-1}\mathbb{I}$.
We additionally include box constraints on the controls: $u(t)\in [-1,1]^2$.
The initial state is $x_0=(-3,-3,0,0,0)$.
The problem is discretized with $\Nfe = 4$, $\ns=4$, $N_s = 30$, and $T=20$.
The full solution to the problem including the initialization via simulation takes $53.58$ seconds.
Several frames of the solution are shown in \Cref{fig:spin} and the resulting controls and state trajectory in \Cref{fig:spinning-states}.

\subsection{Collaborative Pushing}
Finally we describe a problem in which two actuated disc pushers must collaborate in order to transport an ellipsoidal slider.
We include pairwise contacts between each pusher and the slider with a friction coefficient of $\mu = 0.5$ but omit the contact between the two sliders.
The state of this system is represented by $x = (p_1,p_2,p_3,\theta_3)\in\R^7$, i.e., the center of mass positions of the two pushers followed by the center of mass position and angle of the slider.
The system initial state is $x_0 = (-3,3,3,3,0,2,0)$.
We introduce stage and terminal terms $L(x,u) = \transp{u}Q_uu + q_\theta\theta_3^2$ and $M(x(T)) = \transp{(x(T)-\bar{x})}Q_T(x(T)-\bar{x})$ with $Q_u = 10^{-1}I$, $q_\theta = 10^2$, $Q_T=\diag(10^{-1}, 10^{-1},10^{-1}, 10^{-1},10^2,10^2,10^3)$, and target configuration $\bar{x} = (-3,6,3,6,0,5,0)$.
The problem is discretized with $\Nfe = 4$, $\ns=2$, $N_s = 30$, and $T=20$.
The full solution to the problem including the initialization via simulation takes $161.598$ seconds.
Several frames of the solution are given in \Cref{fig:drag}.
\section{Conclusions and Future Work}
In this paper we have proposed a practical direction towards solving planar manipulation problems formulated as projected dynamical systems using direct optimal control.
This approach is extended with a quasi-static friction model modeled via additional complementarity constraints on the tangential velocity.
We demonstrate a method for implicitly defining a signed distance function for ellipsoids via a convex optimization problem.
Finally several examples are described and model solutions are provided for two dimensional pushing problems.
Additional examples of solutions can be found at\footnote{\href{https://youtu.be/T4Fp5GaJC3k}{https://youtu.be/T4Fp5GaJC3k}}.
Future work will encompass extension of this framework to other shapes and the extension of the friction formulation to three dimensions.
\bibliographystyle{plainnat}
\bibliography{syscop}

\end{document}